# Mannheim Offsets of the Timelike Ruled Surfaces with Spacelike Rulings in Dual Lorentzian Space $\mathbb{D}_1^3$

**Mehmet ÖNDER[*], H. Hüseyin UĞURLU[**]**

[*]*Celal Bayar University, Faculty of Science and Arts, Department of Mathematics, Muradiye Manisa.* mehmet.onder@bayar.edu.tr

[**]*Gazi University, Faculty of Education, Department of Secondary Education Science and Mathematics Teaching, Mathematics Teaching Programme, Ankara, Turkey*
hugurlu@gazi.edu.tr,

**Abstract**

In this paper, we obtain the characterizations of Mannheim offsets of the timelike ruled surface with spacelike rulings in dual Lorentzian space $ID_1^3$. We give the relations between terms of their integral invariants and also we give the new characterization of the Mannheim offsets of developable timelike ruled surface. Moreover, we obtain the relationships between the area of projections of spherical images for Mannheim offsets of timelike ruled surfaces and their integral invariants.

## 1. Introduction

Ruled surfaces are the surfaces generated by a continuously moving of a straight line in the space. These surfaces are one of the most important topics of differential geometry. Because of their simple generation, these surfaces arise in a variety of applications including Computer Aided Geometric Design(CAGD), mathematical physics, moving geometry, kinematics for modeling the problems and model-based manufacturing of mechanical products. Especially, the offsets of the surfaces have an important role and applications in these sciences. Some studies dealing with offsets of surfaces have been given in[6,23,24]. Ravani and Ku have defined and given a generalization of the theory of Bertrand curves for Bertrand trajectory ruled surfaces on the line geometry[25]. Küçük and Gürsoy have considered the Bertrand trajectory ruled surfaces in dual space[14]. By using the integral invariants of the ruled surfaces, given in [7,8,9,10], they have given the relations between the invariants of the Bertrand trajectory ruled surfaces.

Moreover, analogue to the notion of Bertrand curve, recently, a new definition of special curves has been given by Liu and Wang: Let $C$ and $C^*$ be two space curves. $C$ is said to be a Mannheim partner curve of $C^*$ if there exists a one to one correspondence between their points such that the binormal vector of $C$ is the principal normal vector of $C^*$ [15]. By considering the study of Ravani and Ku, Orbay and et.al. have given a generalization of the theory of Mannheim curves for ruled surfaces and called Mannheim offsets[16]. Moreover, Önder and et al have studied the Mannheim offsets of the ruled surfaces in dual space and have obtained the relations between the integral invariants of the trajectory Mannheim ruled surfaces[20].

Mannheim offsets of spacelike and timelike ruled surfaces in Minkowski 3-space $E_1^3$ have been studied in [18,19]. In these papers, the authors have given the characterizations of Mannheim offsets of the ruled surfaces in real Lorentzian space(Minkowski 3-space $E_1^3$). They haven't studied the Mannheim offsets of closed trajectory timelike and spacelike ruled surfaces and not considered the dual forms of timelike and spacelike ruled surface.

In this paper, using the dual representations of the timelike ruled surfaces, we give the result characterizing the developable Mannheim offsets of timelike ruled surfaces obtained in [19] in short form. Furthermore, differently from the results given in [19], we examine the Mannheim offsets of trajectory timelike ruled surfaces with spacelike ruling in view of their integral invariants. We obtain that the striction lines of Mannheim offsets of developable timelike trajectory ruled surfaces are Mannheim partner curves in Minkowski 3-space $E_1^3$. Moreover, we give the relations between the integral invariants (such as the angle of pitch and the pitch) of closed trajectory ruled surfaces. Finally, we obtain the relationship between the area of projections of spherical images and integral invariants.

## 2. Differential Geometry of the Ruled Surfaces in Minkowski 3-space

The Minkowski 3-space $IR_1^3$ is the real vector space $IR^3$ provided with the standard flat metric given by
$$\langle \vec{a}, \vec{a} \rangle = -a_1 b_1 + a_2 b_2 + a_3 b_3,$$
where $\vec{a} = (a_1, a_2, a_3)$ and $\vec{b} = (b_1, b_2, b_3) \in IR_1^3$. An arbitrary vector $\vec{a} = (a_1, a_2, a_3)$ of $IR_1^3$ has one of three Lorentzian causal characters; it can be spacelike if $\langle \vec{a}, \vec{a} \rangle > 0$ or $\vec{a} = 0$, timelike if $\langle \vec{a}, \vec{a} \rangle < 0$ and null (lightlike) if $\langle \vec{a}, \vec{a} \rangle = 0$ and $\vec{a} \neq 0$ [17]. Similarly, an arbitrary curve $\vec{\alpha} = \vec{\alpha}(s)$ can locally be spacelike, timelike or null (lightlike), if all of its velocity vectors $\alpha'(s)$ are respectively spacelike, timelike or null (lightlike). We say that a timelike vector is future pointing or past pointing if the first compound of the vector is positive or negative, respectively. The norm of a vector $\vec{a}$ is defined by $\|\vec{a}\| = \sqrt{|\langle \vec{a}, \vec{a} \rangle|}$. Now, let $\vec{a} = (a_1, a_2, a_3)$ and $\vec{b} = (b_1, b_2, b_3)$ be two vectors in $IR_1^3$, then the Lorentzian cross product of $\vec{a}$ and $\vec{b}$ is given by
$$\vec{a} \times \vec{b} = (a_2 b_3 - a_3 b_2, a_1 b_3 - a_3 b_1, a_2 b_1 - a_1 b_2).$$

The sets of the unit timelike and spacelike vectors are called hyperbolic unit sphere and Lorentzian unit sphere, respectively, and denoted by
$$H_0^2 = \{\vec{a} = (a_1, a_2, a_3) \in E_1^3 : \langle \vec{a}, \vec{a} \rangle = -1\}$$
and
$$S_1^2 = \{\vec{a} = (a_1, a_2, a_3) \in E_1^3 : \langle \vec{a}, \vec{a} \rangle = 1\}$$
respectively(See [18,19]).

A surface in the Minkowski 3-space $IR_1^3$ is called a timelike surface if the induced metric on the surface is a Lorentz metric and is called a spacelike surface if the induced metric on the surface is a positive definite Riemannian metric, i.e., the normal vector on the spacelike (timelike) surface is a timelike (spacelike) vector[11].

**Lemma 2.1.** *In the Minkowski 3-space $IR_1^3$, the following properties are satisfied:*
*(i) Two timelike vectors are never orthogonal.*
*(ii) Two null vectors are orthogonal if and only if they are linearly dependent.*
*(iii) A timelike vector is never orthogonal to a null (lightlike) vector* [1].

Let $I$ be an open interval in the real line $IR$, $\vec{k} = \vec{k}(s)$ be a curve in $IR_1^3$ defined on $I$ and $\vec{q} = \vec{q}(s)$ be a unit direction vector of an oriented line in $IR_1^3$. Then we have the following parametrization for a ruled surface $M$,
$$\varphi(s, v) = \vec{k}(s) + v \vec{q}(s). \tag{1}$$

The parametric $s$-curve of this surface is a straight line of the surface which is called ruling. For $v = 0$, the parametric $v$-curve of this surface is $\vec{k} = \vec{k}(s)$ which is called base curve or generating curve of the surface. In particular, if $\vec{q}$ is constant, the ruled surface is said to be cylindrical, and non-cylindrical otherwise.

The striction point on a ruled surface $M$ is the foot of the common normal between two consecutive rulings. The set of the striction points constitutes a curve $\vec{c} = \vec{c}(s)$ lying on the ruled surface and is called striction curve. The parametrization of the striction curve $\vec{c} = \vec{c}(s)$ on a ruled surface is given by

$$\vec{c}(s) = \vec{k}(s) - \frac{\langle d\vec{q}/ds, d\vec{k}/ds \rangle}{\langle d\vec{q}/ds, d\vec{q}/ds \rangle} \vec{q}. \tag{2}$$

So that, the base curve of the ruled surface is its striction curve if and only if $\langle d\vec{q}/ds, d\vec{k}/ds \rangle = 0$.

The distribution parameter (or drall) of the ruled surface in (1) is given by

$$\delta_\varphi = \frac{\left| d\vec{k}/ds, \vec{q}, d\vec{q}/ds \right|}{\langle d\vec{q}/ds, d\vec{q}/ds \rangle} \tag{3}$$

If $\left| d\vec{k}/ds, \vec{q}, d\vec{q}/ds \right| = 0$, then the normal vectors are collinear at all points of the same ruling and at the nonsingular points of the surface $M$, the tangent planes are identical. We then say that the tangent plane contacts the surface along a ruling. Such a ruling is called a *torsal ruling*. If $\left| d\vec{k}/ds, \vec{q}, d\vec{q}/ds \right| \neq 0$, then the tangent planes of the surface $M$ are distinct at all points of the same ruling which is called nontorsal.

In Minkowski 3-space, a ruled surface whose all rulings are torsal is called a *developable ruled surface*. The remaining ruled surfaces are called skew ruled surfaces. Furthermore, the generator $\vec{q}$ of a developable ruled surface is tangent of its striction curve(see [1,18,19] for details).

**Theorem 2.1.** *A ruled surface is developable if and only if at all points the distribution parameter $\delta_\varphi = 0$* [13,26,27].

For the unit normal vector $\vec{m}$ of the ruled surface $M$ we have $\vec{m} = \frac{\vec{\varphi}_s \times \vec{\varphi}_v}{\|\vec{\varphi}_s \times \vec{\varphi}_v\|}$.

So, at the points of a nontorsal ruling $s = s_1$ we have

$$\vec{a} = \lim_{v \to \infty} \vec{m}(s_1, v) = \frac{(d\vec{q}/ds) \times \vec{q}}{\|d\vec{q}/ds\|}.$$

The plane of the ruled surface $M$ which passes through its ruling $s_1$ and is perpendicular to the vector $\vec{a}$ is called the *asymptotic plane* $\alpha$. The tangent plane $\gamma$ passing through the ruling $s_1$ which is perpendicular to the asymptotic plane $\alpha$ is called the *central plane*. Its point of contact $C$ is *central point* of the ruling. The straight lines which pass through point $C$ and are perpendicular to the planes $\alpha$ and $\gamma$ are called the *central tangent* and *central normal*, respectively.

Using the perpendicularly of the vectors $\vec{q}, d\vec{q}/ds$ and the vector $\vec{a}$, the representation of the unit vector $\vec{h}$ of the central normal is given by

$$\vec{h} = \frac{d\vec{q}/ds}{\|d\vec{q}/ds\|}.$$

The orthonormal system $\{C; \vec{q}, \vec{h}, \vec{a}\}$ is called Frenet frame of the ruled surface $M$ such that $\vec{h} = \dfrac{d\vec{q}/ds}{\|d\vec{q}/ds\|}$ and $\vec{a} = \vec{h} \times \vec{q}$ are the central normal and the asymptotic normal direction of $M$, respectively, and $C$ is the striction point.

In [12], Kim and Yoon have given the classifications of the ruled surfaces in Minkowski 3-space. They have defined five types of the ruled surfaces. In these classifications, they have also considered the ruled surfaces with null frenet frames and their null derivatives. Since the angle between null vectors and other vectors is not defined, in this study we don't consider these types of ruled surfaces.

Let now consider the ruled surface $M$ with non-null frenet vectors and their non-null derivatives. According to the Lorentzian characters of the ruling and the central normal vector, we can give the following classifications of the ruled surface $M$ as follows;

**i)** If the central normal vector $\vec{h}$ is spacelike and $\vec{q}$ is timelike, then the ruled surface $M$ is said to be of type $M_-^1$.

**ii)** If the central normal vector $\vec{h}$ and the ruling $\vec{q}$ are both spacelike, then the ruled surface $M$ is said to be of type $M_+^1$.

**iii)** If the central normal vector $\vec{h}$ is timelike and the ruling $\vec{q}$ is spacelike, then the ruled surface $M$ is said to be of type $M_+^2$ [18,19].

The ruled surfaces of type $M_+^1$ and $M_-^1$ are clearly timelike and the ruled surface of type $M_+^2$ is spacelike.

By using these classifications, the parametrization of the ruled surface $M$ can be given as follows,

$$\varphi(s,v) = \vec{k}(s) + v\vec{q}(s), \qquad (4)$$

where $\langle \vec{h}, \vec{h} \rangle = \varepsilon_1 (=\pm 1)$, $\langle \vec{q}, \vec{q} \rangle = \varepsilon_2 (=\pm 1)$.

The set of all bound vectors $\vec{q}(s)$ at the point O constitutes the *directing cone* of the ruled surface $M$. If $\varepsilon_2 = -1$ (resp. $\varepsilon_2 = 1$), the end points of the vectors $\vec{q}(s)$ drive a spherical spacelike (resp. spacelike or timelike) curve on hyperbolic unit sphere $H_0^2$ (resp. on Lorentzian unit sphere $S_1^2$), called the *hyperbolic (resp. Lorentzian) spherical image* of the ruled surface $M$.

Let $\left\{ \vec{q}, \vec{h} = \dfrac{d\vec{q}/ds}{\|d\vec{q}/ds\|}, \vec{a} = \vec{q} \times \vec{h} \right\}$ be a moving othonormal trihedron making a spatial motion along a closed space curve $\vec{k}(s)$, $s \in IR$, in $IR_1^3$ where $\vec{q}$ and $\vec{h}$ are spacelike vectors. In this motion, the oriented line $\vec{q}$ generates a closed timelike ruled surface of the type $M_+^1$ called closed timelike trajectory ruled surface (CTTRS) in $IR_1^3$. A parametric equation of a closed trajectory timelike ruled surface generated by $\vec{q}$-axis is

$$\varphi_q(s,v) = \vec{k}(s) + v\vec{q}(s), \quad \varphi(s+2\pi, v) = \varphi(s,v), \quad s,v \in IR. \qquad (5)$$

Consider the moving orthonormal system $\{\vec{q}, \vec{h}, \vec{a}\}$ which represents a timelike ruled surface of the type $M_+^1$ generated by the vector $\vec{q}$. Then, the axes of the trihedron intersect at the striction point of $\vec{q}$-generator of $\varphi_q$-CTTRS. The structral equations of this motion are

$$\begin{cases} d\vec{q} = k_1 \vec{h} \\ d\vec{h} = -k_1 \vec{q} + k_2 \vec{a} \\ d\vec{a} = k_2 \vec{h} \end{cases} \tag{6}$$

and

$$\frac{db}{ds} = \cosh \sigma \vec{q} + \sinh \sigma \vec{a}, \tag{7}$$

where $b = b(s)$ is the striction line of $\varphi_q$-CTTRS and the differential forms $k_1$, $k_2$ and $\sigma$ are the natural curvature, the natural torsion and the striction of $\varphi_q$-CTTRS, respectively. Here, $s$ is the length of the striction line.

The pole vector and the Steiner vector of the motion are given by

$$\vec{p} = \frac{\vec{\psi}}{\|\vec{\psi}\|}, \quad \vec{d} = \oint \vec{\psi}, \tag{8}$$

respectively, where $\vec{\psi} = k_2 \vec{q} - k_1 \vec{a}$ is the instantaneous Pfaffian vector of the motion[12,19].

The pitch of $\varphi_q$-CTTRS is defined by

$$\ell_q = \oint dv = -\oint \langle d\vec{k}, \vec{q} \rangle. \tag{9}$$

The angle of pitch of $\varphi_q$-CTTRS is given one of the followings

$$\lambda_q = \oint d\theta = \oint \langle d\vec{h}, \vec{a} \rangle = \varepsilon_2 \langle \vec{q}, \vec{d} \rangle = 2\pi - a_q, \tag{10}$$

where $a_q$ is the measure of the spherical surface area bounded by the spherical image of $\varphi_q$-CTTRS. The pitch and the angle of pitch are well-known real integral invariants of closed timelike trajectory ruled surface(See [2,3,4,21,30,31]).

The area vector of a $x$-closed space curve in $IR_1^3$ is given by

$$v_x = \oint x \times dx, \tag{11}$$

and the area of projection of a $x$-closed space curve in direction of the generator of a $y$-CTRS is

$$2f_{x,y} = \langle v_x, y \rangle. \tag{12}$$

(See [14,30]).

## 3. Dual Lorentzian Vectors and E. Study Mapping

In this section, we give a brief summary of the theory of dual numbers and dual Lorentzian vectors.

A dual number has the form $\bar{\lambda} = \lambda + \varepsilon \lambda^*$, where $\lambda$ and $\lambda^*$ are real numbers and $\varepsilon$ is dual unit with $\varepsilon \neq 0$, $\varepsilon^2 = 0$. In differential geometry and motion analysis of spatial mechanisms, a dual number is also commonly referred as a dual angle $\bar{\lambda} = \lambda + \varepsilon \lambda^*$ between two lines in the space. The real part $\lambda$ of the dual angle is the projected angle between the lines, and the dual part $\lambda^*$ is the length along the common normal of the lines.

We denote the set of all dual numbers by $\mathbb{D}$:

$$\mathbb{D} = \{\hat{\lambda} = \lambda + \varepsilon \lambda^* : \lambda, \lambda^* \in IR, \varepsilon^2 = 0\}.$$

Equality, addition and multiplication are defined in $\mathbb{D}$ by

$$\lambda + \varepsilon \lambda^* = \beta + \varepsilon \beta^* \text{ if and only if } \lambda = \beta \text{ and } \lambda^* = \beta^*,$$

$$(\lambda + \varepsilon \lambda^*) + (\beta + \varepsilon \beta^*) = (\lambda + \beta) + \varepsilon (\lambda^* + \beta^*),$$

and

$$(\lambda+\varepsilon\lambda^*)(\beta+\varepsilon\beta^*)=\lambda\beta+\varepsilon(\lambda\beta^*+\lambda^*\beta),$$

respectively. Then it is easy to show that $(\mathbb{D},+,.)$ is a commutative ring with unity. The numbers $\varepsilon\lambda^*$ ($\lambda\in IR$) are divisors of 0.

Now let $f$ be a differentiable function with dual variable $\bar{x}=x+\varepsilon x^*$. Then the Maclaurine series generated by $f$ is

$$f(\bar{x})=f(x+\varepsilon x^*)=f(x)+\varepsilon x^* f'(x),$$

where $f'(x)$ is the derivative of $f$.

Now let $\mathbb{D}^3$ be the set of all triples of dual numbers, i.e.

$$\mathbb{D}^3=\{\tilde{a}=(\bar{a}_1,\bar{a}_2,\bar{a}_3)\mid \bar{a}_i\in\mathbb{D},\ 1\leq i\leq 3\}.$$

The elements of $\mathbb{D}^3$ are called as dual vectors. A dual vector $\tilde{a}$ may be expressed in the form $\tilde{a}=\vec{a}+\varepsilon\vec{a}^*$, where $a$ and $a^*$ are the vectors of $IR^3$.

Now let $\tilde{a}=\vec{a}+\varepsilon\vec{a}^*$, $\tilde{b}=\vec{b}+\varepsilon\vec{b}^*\in\mathbb{D}^3$ and $\bar{\lambda}=\lambda_1+\varepsilon\lambda_1^*\in\mathbb{D}$. Then we define

$$\tilde{a}+\tilde{b}=\vec{a}+\vec{b}+\varepsilon(\vec{a}^*+\vec{b}^*),$$
$$\bar{\lambda}\tilde{a}=\lambda_1\vec{a}+\varepsilon(\lambda_1\vec{a}^*+\lambda_1^*\vec{a}).$$

Then, $\mathbb{D}^3$ becomes a unitary module with these operations. It is called $\mathbb{D}$-module or dual space[5,10,32].

The Lorentzian inner product of two dual vectors $\tilde{a}=\vec{a}+\varepsilon\vec{a}^*$, $\tilde{b}=\vec{b}+\varepsilon\vec{b}^*\in\mathbb{D}^3$ is defined by

$$\langle\tilde{a},\tilde{b}\rangle=\langle\vec{a},\vec{b}\rangle+\varepsilon\left(\langle\vec{a},\vec{b}^*\rangle+\langle\vec{a}^*,\vec{b}\rangle\right),$$

where $\langle\vec{a},\vec{b}\rangle$ is the Lorentzian inner product of the vectors $\vec{a}$ and $\vec{b}$ in the Minkowski 3-space $IR_1^3$. Then a dual vector $\tilde{a}=\vec{a}+\varepsilon\vec{a}^*$ is said to be dual timelike if $\vec{a}$ is timelike, dual spacelike if $\vec{a}$ is spacelike or $\vec{a}=0$ and dual lightlike (null) if $\vec{a}$ is lightlike (null) and $\vec{a}\neq 0$.

The set of all dual Lorentzian vectors is called dual Lorentzian space and it is denoted by $\mathbb{D}_1^3$:

$$\mathbb{D}_1^3=\{\tilde{a}=\vec{a}+\varepsilon\vec{a}^*:\vec{a},\vec{a}^*\in IR_1^3\}.$$

The Lorentzian cross product of dual vectors $\tilde{a}$ and $\tilde{b}\in\mathbb{D}_1^3$ is defined by

$$\tilde{a}\times\tilde{b}=\vec{a}\times\vec{b}+\varepsilon(\vec{a}^*\times\vec{b}+\vec{a}\times\vec{b}^*),$$

where $\vec{a}\times\vec{b}$ is the Lorentzian cross product in $IR_1^3$.

Let $\tilde{a}=\vec{a}+\varepsilon\vec{a}^*\in\mathbb{D}_1^3$. Then $\tilde{a}$ is said to be dual unit timelike(resp. spacelike) vector if the vectors $\vec{a}$ and $\vec{a}^*$ satisfy the following equations:

$<\vec{a},\vec{a}>=-1$ (resp. $<\vec{a},\vec{a}>=1$), $<\vec{a},\vec{a}^*>=0.$

The set of all unit dual timelike vectors is called the dual hyperbolic unit sphere, and is denoted by $\tilde{H}_0^2$. Similarly, the set of all unit dual spacelike vectors is called the dual Lorentzian unit sphere, and is denoted by $\tilde{S}_1^2$ [30,31,33].

***Theorem 3.1. (E. Study's Mapping):*** *The dual timelike(respectively spacelike) unit vectors of the dual hyperbolic(respectively Lorentzian) unit sphere $\tilde{H}_0^2$ (respectively $\tilde{S}_1^2$) are in on-to-one correspondence with the directed timelike(respectively spacelike) lines of the Minkowski 3-space $IR_1^3$ [33].*

**Definition 3.1. i) Dual Hyperbolic angle:** Let $\tilde{x}$ and $\tilde{y}$ be dual timelike vectors in $ID_1^3$. Then the dual angle between $\tilde{x}$ and $\tilde{y}$ is defined by $<\tilde{x},\tilde{y}>= -\|\tilde{x}\|\|\tilde{y}\|\cosh\bar{\theta}$. The dual number $\bar{\theta} = \theta + \varepsilon\theta^*$ is called the dual *hyperbolic angle*.

**ii) Dual Central angle:** Let $\tilde{x}$ and $\tilde{y}$ be dual spacelike vectors in $ID_1^3$ that span a dual timelike vector subspace. The dual angle between $\tilde{x}$ and $\tilde{y}$ is defined by $<\tilde{x},\tilde{y}>= \|\tilde{x}\|\|\tilde{y}\|\cosh\bar{\theta}$. The dual number $\bar{\theta} = \theta + \varepsilon\theta^*$ is called the dual *central angle*.

**iii) Dual Spacelike angle:** Let $\tilde{x}$ and $\tilde{y}$ be dual spacelike vectors in $ID_1^3$ that span a dual spacelike vector subspace. Then the angle between $\tilde{x}$ and $\tilde{y}$ is defined by $<\tilde{x},\tilde{y}>= \|\tilde{x}\|\|\tilde{y}\|\cos\bar{\theta}$. The dual number $\bar{\theta} = \theta + \varepsilon\theta^*$ is called the dual *spacelike angle*.

**iv) Dual Lorentzian timelike angle:** Let $\tilde{x}$ be a dual spacelike vector and $\tilde{y}$ be a dual timelike vector in $ID_1^3$. Then the angle between $\tilde{x}$ and $\tilde{y}$ is defined by $<\tilde{x},\tilde{y}>= \|\tilde{x}\|\|\tilde{y}\|\sinh\bar{\theta}$. The dual number $\bar{\theta} = \theta + \varepsilon\theta^*$ is called the dual *Lorentzian timelike angle*[29,30,31].

Let $\tilde{K}$ be a moving dual Lorentzian unit sphere generated by a dual orthonormal system

$$\left\{\tilde{q}, \tilde{h} = \frac{d\tilde{q}}{\|d\tilde{q}\|}, \tilde{a} = \tilde{q}\times\tilde{h}\right\}, \tilde{q} = \vec{q} + \varepsilon\vec{q}^*, \tilde{h} = \vec{h} + \varepsilon\vec{h}^*, \tilde{a} = \vec{a} + \varepsilon\vec{a}^*, \qquad (13)$$

and $\tilde{K}'$ be a fixed dual Lorentzian unit sphere with the same center where $\tilde{q}, \tilde{h}$ are spacelike vectors. Then, the derivative equations of the dual spherical closed motion of $\tilde{K}$ with respect to $\tilde{K}'$ are

$$\begin{cases} d\tilde{q} = \bar{k}_1\tilde{h} \\ d\tilde{h} = -\bar{k}_1\tilde{q} + \bar{k}_2\tilde{a} \\ d\tilde{a} = \bar{k}_2\tilde{h} \end{cases} \qquad (14)$$

where $\bar{k}_1(s) = k_1(s) + \varepsilon k_1^*(s)$, $\bar{k}_2(s) = k_2(s) + \varepsilon k_2^*(s)$, $(s \in IR)$ are dual curvature and dual torsion, respectively. From the E. Study mapping, during the spherical motion of $\tilde{K}$ with respect to $\tilde{K}'$, the dual unit vector $\tilde{q}$ draws a dual curve on dual unit Lorentzian sphere $K'$ and this curve represents a timelike ruled surface with spacelike ruling $\vec{q}$ in line space $IR_1^3$.

Dual vector $\tilde{\psi} = \vec{\psi} + \varepsilon\vec{\psi}^* = \bar{k}_2\tilde{q} - \bar{k}_1\tilde{a}$ is called the instantaneous Pfaffian vector of the motion and the vector $\tilde{P}$ given by $\tilde{\psi} = \|\tilde{\psi}\|\tilde{P}$ is called the dual pole vector of the motion. Then the vector

$$\tilde{d} = \oint \tilde{\psi}, \qquad (15)$$

is called the dual Steiner vector of the closed motion.

By considering the E. Study mapping, the dual equations (14) correspond to the real equations (6) and (7) of a closed spatial motion in $IR_1^3$. So, the differentiable dual closed curve $\tilde{q} = \tilde{q}(s)$ is corresponds to a closed timelike trajectory ruled surface with spacelike ruling in the line space $IR_1^3$ and denoted by $\varphi_q$-CTTRS.

A dual integral invariant of a $\varphi_q$-CTTRS can be given in terms of real integral invariants as follows and is called the dual angle of pitch of a $\varphi_q$-CTTRS

$$\bar{\wedge}_q = \oint \langle d\tilde{h}, \tilde{a}\rangle = -\langle \tilde{q}, \tilde{d}\rangle = 2\pi - \bar{a}_q = \lambda_q + \varepsilon\ell_q \qquad (16)$$

where $\tilde{d} = \vec{d} + \varepsilon \vec{d}^*$ and $\bar{a}_q = a_q + \varepsilon a_q^*$ are the dual Steiner vector of the motion and the measure of dual spherical surface area of $\varphi_q$-CTTRS, respectively. Here, $\varepsilon$ is dual unit.

Analogue to the real area vector given in (11), the dual area vector of a $\tilde{q}$-closed dual curve is given by

$$\tilde{w}_{\tilde{q}} = \oint \tilde{q} \times d\tilde{q} \qquad (17)$$

and the dual area of projection of a $\tilde{q}$-closed dual curve in direction of the generator of a $\tilde{q}_1$-CTRS is

$$2\bar{f}_{\tilde{q},\tilde{q}_1} = \langle \tilde{w}_{\tilde{q}}, \tilde{q}_1 \rangle. \qquad (18)$$

(See [21,30,31] for details).

**4. Mannheim Offsets of Timelike Trajectory Ruled Surfaces with Spacelike Rulings**

Let $\varphi_q$ be a timelike trajectory ruled surfaces of the type $M_+^1$ generated by dual spacelike vectors $\tilde{q}$ and let the dual orthonormal frame of $\varphi_q$ be $\{\tilde{q}(s), \tilde{h}(s), \tilde{a}(s)\}$. The trajectory ruled surface $\varphi_{q_1}$, generated by dual vector $\tilde{q}_1$, with dual orthonormal frame $\{\tilde{q}_1(s_1), \tilde{h}_1(s_1), \tilde{a}_1(s_1)\}$ is said to be Mannheim offset of the timelike trajectory ruled surface $\varphi_q$, if

$$\tilde{a}(s) = \tilde{h}_1(s_1) \qquad (19)$$

where $s$ and $s_1$ are the arc-length of the striction lines of $\varphi_q$ and $\varphi_{q_1}$, respectively. By this definition and considering the classifications of the ruled surfaces in Minkowski 3-space, we have that the Mannheim offsets $\varphi_{q_1}$ is a spacelike ruled surface of the type $M_+^2$. Thus by the Definition 3.1, we can write

$$\begin{pmatrix} \tilde{q}_1 \\ \tilde{h}_1 \\ \tilde{a}_1 \end{pmatrix} = \begin{pmatrix} \cos\bar{\theta} & \sin\bar{\theta} & 0 \\ 0 & 0 & 1 \\ \sin\bar{\theta} & -\cos\bar{\theta} & 0 \end{pmatrix} \begin{pmatrix} \tilde{q} \\ \tilde{h} \\ \tilde{a} \end{pmatrix}. \qquad (20)$$

In (20), $\bar{\theta} = \theta + \varepsilon\theta^*$, $(\theta, \theta^* \in IR)$ is the dual angle between the generators $\tilde{q}$ and $\tilde{q}_1$ of the Mannheim trajectory ruled surface $\varphi_q$ and $\varphi_{q_1}$. The angle $\theta$ is called the offset angle and $\theta^*$ is called the offset distance. Then, $\bar{\theta} = \theta + \varepsilon\theta^*$ is called dual offset angle of the Mannheim trajectory ruled surface $\varphi_q$ and $\varphi_{q_1}$. If $\theta = 0$ and $\theta = \pi/2$, then the Mannheim offsets are said to be oriented offsets and right offsets, respectively. Thus, we can give the followings.

In the latter of the paper, by $\varphi_q$ and $\varphi_{q_1}$, we will mean a timelike ruled surface of the type $M_+^1$ and a spacelike ruled surface of the type $M_+^2$, respectively.

**Theorem 4.1.** *Let $\varphi_q$ and $\varphi_{q_1}$ be the Mannheim trajectory ruled surfaces. The offset angle and the offset distance are given by*

$$\theta = -\int k_1 ds, \quad \theta^* = -\int k_1^* ds \qquad (21)$$

*respectively.*

**Proof:** Let $\varphi_q$ and $\varphi_{q_1}$ form a Mannheim offset. From (20) we have

$$\tilde{q}_1 = \cos\bar{\theta}\tilde{q} + \sin\bar{\theta}\tilde{h}. \qquad (22)$$

Differentiating (22) and by using (14) and (19), we may write

$$\frac{d\tilde{q}_1}{ds} = -\left(\frac{d\bar{\theta}}{ds} + \bar{k}_1\right)\tilde{a}_1 + \bar{k}_2 \sin\bar{\theta}\tilde{h}_1. \qquad (23)$$

Since $\dfrac{d\tilde{q}_1}{ds}$ is orthogonal to $\tilde{a}_1$, from (23) we get
$$\bar{\theta} = -\int \bar{k}_1 ds.$$
Separating the last equation into real and dual parts we have
$$\theta = -\int k_1 ds, \quad \theta^* = -\int k_1^* ds. \tag{24}$$

**Theorem 4.2.** *The closed timelike trajectory ruled surface $\varphi_q$ and the closed spacelike trajectory ruled surface $\varphi_{q_1}$ form a Mannheim offsets if and only if the following relationship holds*
$$\overline{\wedge}_{q_1} = \overline{\wedge}_q \cos\bar{\theta} + \overline{\wedge}_h \sin\bar{\theta}, \quad \bar{\theta} = -\int \bar{k}_1 ds. \tag{25}$$

**Proof:** Let the $\varphi_q$-CTTRS and the $\varphi_{q_1}$-CSTRS form a Mannheim offsets. Then, by direct calculation, from (16) and (20), the dual angle of pitch of $\varphi_{q_1}$-CSTRS is given by
$$\overline{\wedge}_{q_1} = \overline{\wedge}_q \cos\bar{\theta} + \overline{\wedge}_h \sin\bar{\theta}, \quad \bar{\theta} = -\int \bar{k}_1 ds.$$

Conversely, if (25) holds, it is easily seen that $\varphi_q$ and $\varphi_{q_1}$-CTRS form a Mannheim offsets.

Equality (25) is a dual characterization of Mannheim offsets of CTRS in terms of their dual integral invariants. Separating (25) into real and dual parts, we obtain
$$\begin{cases} \lambda_{q_1} = \lambda_q \cos\theta + \lambda_h \sin\theta \\ \ell_{q_1} = (\ell_q + \theta^* \lambda_h)\cos\theta - (\ell_h + \theta^* \lambda_q)\sin\theta \end{cases} \tag{26}$$
From (26), we have the following special cases:

**Case 1.** If $\varphi_q$ and $\varphi_{q_1}$ are the oriented closed Mannheim trajectory ruled surfaces i.e., $\theta = 0$, then the relationships between the real integral invariants of $\varphi_q$ and $\varphi_{q_1}$-CTRS are given as follows,
$$\lambda_{q_1} = \lambda_q, \quad \ell_{q_1} = \ell_q + \theta^* \lambda_h. \tag{27}$$
Furthermore, the measure of spherical surface areas bounded by the spherical images of $\varphi_q$ and $\varphi_{q_1}$-CTRS Mannheim offsets and $\varphi_h$-CTRS are given by
$$a_{q_1} = a_q \text{ and } a_{q_1}^* = -a_q^* + \theta^*(2\pi - a_h) \tag{28}$$

**Case 2.** If $\varphi_q$ and $\varphi_{q_1}$ are the right closed Mannheim trajectory ruled surfaces i.e., $\theta = \pi/2$, then the relationships between the real integral invariants of $\varphi_q$ and $\varphi_{q_1}$-CTRS are given as follows
$$\lambda_{q_1} = \lambda_h, \quad \ell_{q_1} = -\ell_h - \theta^* \lambda_q. \tag{29}$$
Then, the measure of spherical surface areas bounded by the spherical images of $\varphi_{q_1}$, $\varphi_{q_1}$ and $\varphi_h$-CTRS's are
$$a_{q_1} = a_h, \quad a_{q_1}^* = -a_h^* + (2\pi - a_q)\theta^*. \tag{30}$$

**Case 3.** If $\theta^* = 0$, i.e., the generators $\vec{q}$ and $\vec{q}_1$ of the Mannheim offset surfaces intersect, then we have
$$\begin{cases} \lambda_{q_1} = \lambda_q \cos\theta + \lambda_h \sin\theta \\ \ell_{q_1} = \ell_q \cos\theta - \ell_h \sin\theta \end{cases} \tag{31}$$

In this case, $\varphi_q$ and $\varphi_{q_1}$-CTRS's are intersect along their striction lines. It means, their striction lines are the same.

Let now consider that what the condition for the developable Mannheim offset of a CTTRS is. Let the timelike trajectory ruled surface $\varphi_q$ with spacelike ruling $\vec{q}$ and the spacelike trajectory ruled surface $\varphi_{q_1}$ be the Mannheim offset surfaces and let $\vec{\alpha}(s)$ and $\vec{\beta}(s_1)$ be the striction lines of $\varphi_q$ and $\varphi_{q_1}$-CTRS, respectively. Then, we can write

$$\vec{\beta}(s) = \vec{\alpha}(s) + \theta^* \vec{a}(s) \tag{32}$$

where $s$ is the arc-length of $\vec{\alpha}(s)$. Assume that $\varphi_q$-CTTRS is developable. Then from (3) and (7) we have

$$\delta_q = \frac{\langle \cosh \sigma \vec{q} + \sinh \sigma \vec{a},\ \vec{q} \times k_1 \vec{h} \rangle}{\langle k_1 \vec{h}, k_1 \vec{h} \rangle} = -\frac{\sinh \sigma}{k_1} = 0 \tag{33}$$

Then, we have $\sigma = 0$. Thus, from (7)

$$\frac{d\vec{\alpha}}{ds} = \vec{q}. \tag{34}$$

Hence, along the striction line $\vec{\alpha}(s)$, the orthogonal frame $\{\vec{q}, \vec{h}, \vec{a}\}$ coincides with the Frenet frame $\{\vec{T}, \vec{N}, \vec{B}\}$ and the differential forms $k_1$ and $k_2$ turn into the curvature $\kappa_\alpha$ and torsion $\tau_\alpha$ of the striction line $\alpha(s)$, respectively. Then, by the aid of (7), (32) and (34) we have

$$\frac{d\vec{\beta}}{ds} = \vec{q} + \theta^* \tau_\alpha \vec{h}. \tag{35}$$

On the other hand, from (19) and (20) we obtain

$$\frac{d\vec{q}_1}{ds} = -\left(\frac{d\theta}{ds} + \kappa_\alpha\right) \sin \theta \vec{q} + \left(\frac{d\theta}{ds} + \kappa_\alpha\right) \cos \theta \vec{h} + \tau_\alpha \sin \theta \vec{a}. \tag{36}$$

By using (21) and the fact that $k_1 = \kappa_\alpha$, from (36) we have

$$\frac{d\vec{q}_1}{ds} = \tau_\alpha \sin \theta \vec{a} \tag{37}$$

From (35) and (37) we have

$$\delta_{q_1} = \frac{\langle d\vec{\beta},\ \vec{q}_1 \times d\vec{q}_1 \rangle}{\langle d\vec{q}_1, d\vec{q}_1 \rangle} = \frac{\sin \theta - \theta^* \tau_\alpha \cos \theta}{\tau_\alpha \sin \theta}. \tag{38}$$

Thus, from (34) and (38), it can be stated that if the Mannheim offsets of $\varphi_q$ and $\varphi_{q_1}$ are developable then the following relationship holds

$$\sin \theta - \theta^* \tau_\alpha \cos \theta = 0. \tag{39}$$

Equation (39) characterizes the developable Mannheim offsets of timelike ruled surface and it is also given in [12], Theorem 4.2, with another proof.

If (39) holds, along the striction line $\beta(s_1)$, the orthogonal frame $\{\vec{q}_1, \vec{h}_1, \vec{a}_1\}$ coincides with the Frenet frame $\{\vec{T}_1, \vec{N}_1, \vec{B}_1\}$. Thus, the following theorem may be given.

**Theorem 4.3.** *If $\varphi_q$ and $\varphi_{q_1}$ are Mannheim offsets of developable trajectory ruled surfaces then their striction lines are Mannheim partner curves in Minkowski 3-space $E_1^3$.*

From (39) we can give the following special case:

**Case 4.** If $\theta^* = 0$, from (39) we have $\theta = 0$. Then the trajectory ruled surfaces $\varphi_q$ and $\varphi_{q_1}$ coincide.

If $\varphi_q$-CTTRS is developable then from the equations (9), (20) and (35) the pitch $\ell_{q_1}$ of $\varphi_{q_1}$-CSTRS is

$$\ell_{q_1} = -\oint (\cos\theta + \theta^* \tau_\alpha \sin\theta) ds$$

Then, we can give the following:

**Case 5.** If $\varphi_q$-CTTRS is developable then the relation between the pitch $\ell_{q_1}$ of $\varphi_{q_1}$-CSTRS and the torsion of striction line $\alpha(s)$ of $\varphi_q$-CTTRS is given by

$$\ell_{q_1} = -\oint (\cos\theta + \theta^* \tau_\alpha \sin\theta) ds \tag{40}$$

Let now consider the area of projections of Mannheim offsets. From (18), the dual area vectors of the spherical images of $\varphi_q$ and $\varphi_{q_1}$ are

$$\begin{cases} \tilde{w}_q = -\tilde{d} + \overline{\lambda}_q \tilde{q}, \\ \tilde{w}_{q_1} = -\tilde{d} + \overline{\lambda}_{q_1} \tilde{q}_1, \end{cases} \tag{41}$$

respectively. Then, from (18), the dual area of projection of the spherical image of $\varphi_{q_1}$-CSTRS in the direction $\vec{q}$, generators of $\varphi_q$-offsets, is

$$2\overline{f}_{\tilde{q}_1,\tilde{q}} = \langle \tilde{w}_{q_1}, \tilde{q} \rangle = \overline{\lambda}_q + \overline{\lambda}_{q_1} \cos\overline{\theta} \tag{42}$$

Separating (42) into real and dual parts we have the following theorem.

***Theorem 4.4.*** *Let $\varphi_q$-CTTRS and $\varphi_{q_1}$-CSTRS form a Mannheim offset. The relationships between the area of projections of spherical images of the Mannheim offsets $\varphi_q$ and $\varphi_{q_1}$ and their integral invariants are given as follows,*

$$\begin{cases} 2f_{q_1,q} = \lambda_q + \lambda_{q_1} \cos\theta, \\ 2f^*_{q_1,q} = -\ell_q + \ell_{q_1} \cos\theta + \lambda_{q_1} \theta^* \sin\theta, \end{cases} \tag{43}$$

Then we obtain the followings:

**Case 6.** If $\varphi_q$ and $\varphi_{q_1}$-CTRS are the oriented surfaces, i.e., $\theta = 0$, then from (43) we have

$$2f_{q_1,q} = \lambda_q + \lambda_{q_1}, \quad 2f^*_{q_1,q} = \ell_q - \ell_{q_1} \tag{44}$$

**Case 7.** If $\varphi_q$ and $\varphi_{q_1}$ are the right closed Mannheim trajectory ruled surfaces i.e., $\theta = \pi/2$, then from (44) we have

$$2f_{q_1,q} = \lambda_q, \quad 2f^*_{q_1,q} = \ell_q - \lambda_{q_1} \theta^*. \tag{45}$$

Similarly, the dual area of projection of spherical image of $\varphi_{q_1}$-CSTRS in the direction $\tilde{h}$ is

$$2\overline{f}_{\tilde{q}_1,\tilde{h}} = \langle \tilde{w}_{q_1}, \tilde{h} \rangle = \overline{\lambda}_h + \overline{\lambda}_{q_1} \sin\overline{\theta} \tag{46}$$

Separating (46) into real and dual parts we have the following theorem:

**Theorem 4.5.** *Let $\varphi_q$-CTTRS and $\varphi_{q_1}$-CSTRS form a Mannheim offset. The relationships between the area of projections of spherical images of the trajectory ruled surfaces $\varphi_h$ and $\varphi_{q_1}$ and their integral invariants are given as follows*

$$\begin{cases} 2f_{q_1,h} = \lambda_h + \lambda_{q_1} \sin\theta, \\ 2f^*_{q_1,h} = -\ell_h - \ell_{q_1} \sin\theta + \lambda_{q_1}\theta^* \cos\theta, \end{cases} \quad (47)$$

Then we obtain the followings special cases:

**Case 8.** If $\varphi_q$ and $\varphi_{q_1}$-CTRS are the oriented Mannheim surfaces, i.e., $\theta = 0$, then from (47) we have

$$2f_{q_1,h} = \lambda_h, \quad 2f^*_{q_1,h} = -\ell_h + \lambda_{q_1}\theta^* \quad (48)$$

**Case 9.** If $\varphi_q$ and $\varphi_{q_1}$ are the right closed Mannheim trajectory ruled surfaces i.e., $\theta = \pi/2$, then from (47) we have

$$2f_{q_1,h} = \lambda_h + \lambda_{q_1}, \quad 2f^*_{q_1,h} = -\ell_h - \ell_{q_1} \quad (49)$$

The dual area of projection of spherical image of $\varphi_{q_1}$-CSTRS in the direction $\tilde{a}$ is

$$2\overline{f}_{\tilde{q}_1,\tilde{a}} = \langle \tilde{w}_{q_1}, \tilde{a} \rangle = \overline{\wedge}_a = \overline{\wedge}_{h_1} = 0 \quad (50)$$

Separating (50) into real and dual parts we have the following theorem:

**Theorem 4.6.** *Let $\varphi_q$-CTTRS and $\varphi_{q_1}$-CSTRS form a Mannheim offset. The relationships between the area of projections of spherical images of the trajectory ruled surfaces $\varphi_a$ and $\varphi_{q_1}$ and their integral invariants are given as follows*

$$2f_{q_1,a} = \lambda_a = \lambda_{h_1} = 0, \quad 2f^*_{q_1,h} = -\ell_a - \ell_{h_1} = 0. \quad (51)$$

**5. Example:** Let consider the timelike ruled surface with spacelike ruling given by the parametrization

$$\varphi(s,v) = (0, \cos s, \sin s) + v(c, -\sin s, \cos s), \quad |c| < 1 \quad (52)$$

which corresponds to the dual curve

$$\tilde{q} = \frac{1}{\sqrt{1-c^2}}(c, -\sin s, \cos s) + \varepsilon \frac{1}{\sqrt{1-c^2}}(1, c\sin s, c\cos s).$$

(See Fig. 1.). Then, the dual frame of the surface $\varphi(s,v)$ is given as follows,

$$\tilde{q} = \frac{1}{\sqrt{1-c^2}}(c, -\sin s, \cos s) + \varepsilon \frac{1}{\sqrt{1-c^2}}(1, c\sin s, c\cos s)$$

$$\tilde{h} = (0, -\cos s, -\sin s) + \varepsilon 0$$

$$\tilde{a} = \frac{1}{\sqrt{1-c^2}}(1, c\sin s, c\cos s) + \varepsilon \frac{1}{\sqrt{1-c^2}}\left(c(\cos^2 s - \sin^2 s), \sin s, \cos s\right)$$

Some Mannheim offsets of (52) can be given as follows:

**i)** The oriented Mannheim offset of (52) with dual offset angle $\overline{\theta} = 0 + \varepsilon\sqrt{1-c^2}$ is given by

$$\varphi_1(s,v) = (1, \cos s + c\sin s, \sin s + c\cos s) + v\frac{1}{\sqrt{1-c^2}}(c, -\sin s, \cos s) \quad (53)$$

(See Fig. 2).

**ii)** The right Mannheim offset of (52) with dual offset angle $\bar{\theta} = \pi/2 + \varepsilon\sqrt{1-c^2}$ is given by

$$\varphi_2(s,v) = (1, \cos s + c \sin s, \sin s + c \cos s) + v\left(\frac{c}{\sqrt{1-c^2}}, -\cos s, -\sin s\right) \quad (54)$$

which is rendered in Fig. 3.

**iii)** The Mannheim offset of (52) with dual offset angle $\bar{\theta} = \pi/3 + \varepsilon\sqrt{1-c^2}\,s$ is given by the parametrization

$$\varphi_3(s,v) = (s, \cos s + cs \sin s, \sin s + cs \cos s)$$
$$+ v\left(\frac{c}{2\sqrt{1-c^2}}, -\frac{c}{2\sqrt{1-c^2}}\sin s - \frac{\sqrt{3}}{2}\cos s, \frac{1}{2\sqrt{1-c^2}}\cos s - \frac{\sqrt{3}}{2}\sin s\right). \quad (55)$$

which is rendered in Fig. 4.

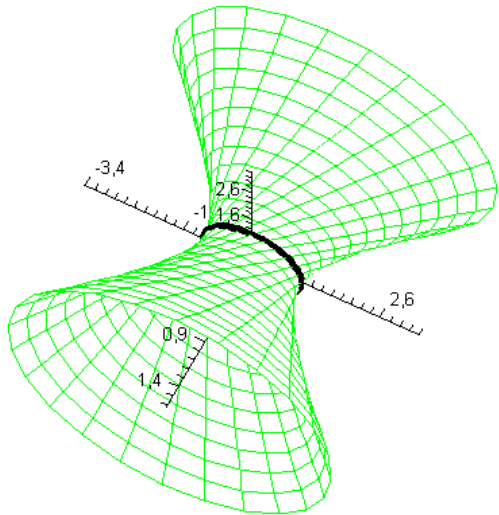

**Fig. 1.** The surface $\varphi(s,v)$.

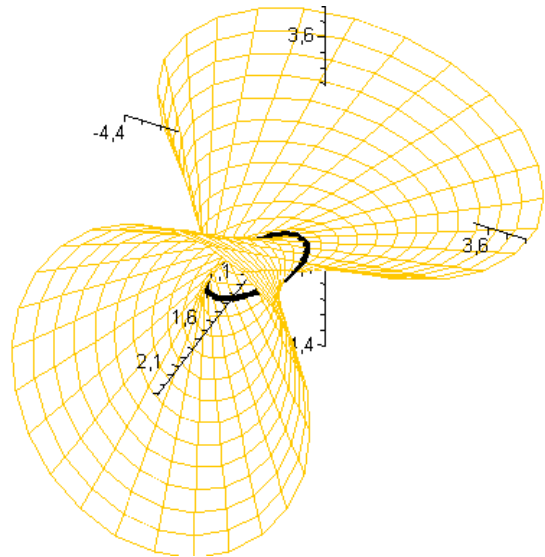

**Fig. 2.** The Mannheim offset $\varphi_1(s,v)$ of $\varphi(s,v)$ with dual offset angle $\bar{\theta} = 0 + \varepsilon\sqrt{1-c^2}$.

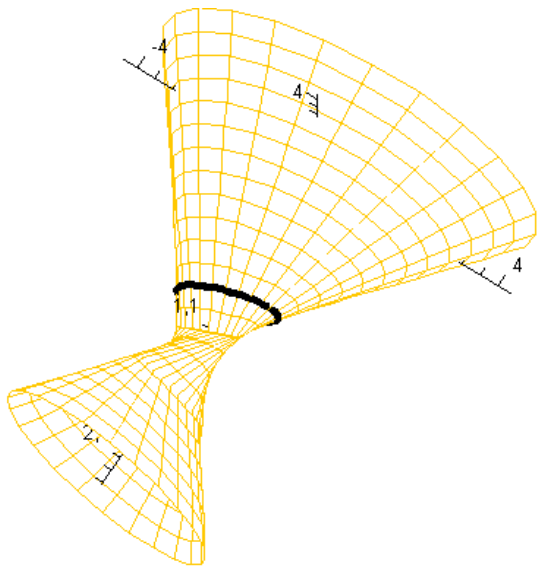

**Fig. 3.** The Mannheim offset $\varphi_2(s,v)$ of $\varphi(s,v)$ with dual offset angle $\bar{\theta} = \pi/2 + \varepsilon\sqrt{1-c^2}$.

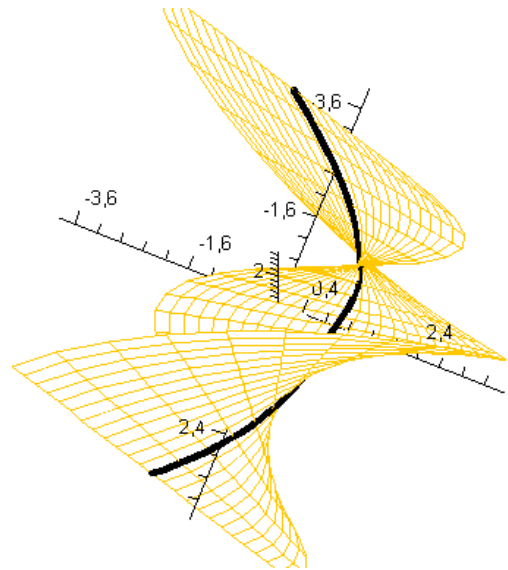

**Fig. 4.** The Mannheim offset $\varphi_3(s,v)$ of $\varphi(s,v)$ with dual offset angle $\bar{\theta} = \pi/3 + \varepsilon\sqrt{1-c^2}\,s$.

In Figures, the curves rendered in black are the striction lines of the offset surfaces.

## 6. Conclusion

In this paper, we give the characterizations of Mannheim offsets of timelike ruled surfaces with spacelike ruling. We obtain the relations between the invariants of Mannheim offsets of timelike ruled surfaces and show that the striction lines of the Mannheim offsets of the developable ruled surface are Mannheim partner curves in Minkowski 3-space $E_1^3$.